%% file: arxiv_RQMC_weak_derivatives.tex
\documentclass[a4paper,11pt]{amsart}
\usepackage[utf8]{inputenc}

\usepackage{enumitem}

\usepackage{a4wide}

\usepackage{orcidlink}
\usepackage{amsmath,amssymb}
\usepackage{bm}
\usepackage{multirow}
\usepackage{algpseudocode}

\usepackage{hyperref}

\usepackage[misc]{ifsym} 

\usepackage{color}
\usepackage{graphicx}
\usepackage{subcaption}

\usepackage{amsfonts}
\usepackage{amsthm}
\usepackage{mathtools}

\usepackage[mathscr]{eucal}

\usepackage{pgfplots}
\pgfplotsset{compat=1.18}

\usepackage{bbm}
\ifpdf
\DeclareGraphicsExtensions{.eps,.pdf,.png,.jpg}
\else
\DeclareGraphicsExtensions{.eps}
\fi

\newtheorem{remark}{Remark}

\input{commands_qmc.tex}

\title{Integrability of weak mixed first-order derivatives and convergence rates of scrambled digital nets}

\begin{document}

\author[Y.~Liu]{Yang Liu}
\address[Y.~Liu]{CEMSE, King Abdullah University of Science and Technology, Thuwal, Saudi Arabia} \email[]{yang.liu.3@kaust.edu.sa}




\maketitle
\begin{abstract}
	We consider the $L^p$ integrability of weak mixed first-order derivatives of the integrand and study convergence rates of scrambled digital nets. We show that the generalized Vitali variation with parameter $\alpha \in [\frac{1}{2}, 1]$ from [Dick and Pillichshammer, 2010] is bounded above by the $L^p$ norm of the weak mixed first-order derivative, where $p = \frac{2}{3-2\alpha}$. Consequently, when the weak mixed first-order derivative belongs to $L^p$ for $1 \leq p \leq 2$, the variance of the scrambled digital nets estimator convergences at a rate of $\mathcal{O}(N^{-4+\frac{2}{p}} \log^{s-1} N)$. Numerical experiments further validate the theoretical results.
\end{abstract}

\section{Introduction}

In this short note, we analyze the variance of the scrambled digital net estimator. Specifically, we consider the $L^p$ integrability of weak, mixed first-order derivatives, where $1\leq p \leq 2$, and establish connections to the generalized Vitali variation proposed in~\cite{dick2010digital}.

Several previous works have studied the convergence rate of the scrambled Sobol' sequence and achieved the rate $\mathcal{O}(N^{-3}\log^{s-1} N)$ for estimator variance under different conditions: Yue and Mao~\cite{yue1999variance} propose a generalized Lipschitz continuity condition on the integrand, whereas Owen~\cite{owen2008local} requires the mixed first-order derivative to be in $L^{\infty}$. Dick and Pillichshammer~\cite{dick2010digital} introduce a generalized Vitali variation to study the scrambled digital net estimator variance and demonstrate that $L^2$ integrability of the mixed first-order derivative is sufficient for achieving the aforementioned convergence rate when the derivative is continuous. In addition, Liu~\cite{liu2024randomized} studies integrands satisfying boundary growth conditions characterized by parameter $A^{*}$ and proves convergence rates of $\mathcal{O}(N^{-2+2A^{*}} \log^{s-1} N)$ for $-1/2 < A^{*} < 1/2$. 

In this work, we establish an estimator variance convergence rate of $\mathcal{O}(N^{-4+\frac{2}{p}} \log^{s-1} N)$ for integrands whose weak derivative $\partial^{1:s} f \in L^{p}$, where $1 \leq p \leq 2$. In the rest of this short note, we introduce necessary notations for the scrambled digital nets and provide convergence analysis, as well as some discussions in Section~\ref{sec:background}. We present some numerical results in~\ref{sec:numex} to support our theoretical results.
\section{Scrambled digital nets and convergence rates}
\label{sec:background}
In this section, we introduce notations for scrambled digital nets and analyze the variance of the estimator.

\subsection{Digital nets}
We begin with the definition of a $(t, m, s)$-net in base $b$. A $(t, m, s)$-net in base $b$ is a set of $b^m$ points in $[0, 1)^s$ such that every $s$-dimensional elementary interval
	\begin{equation}
		\label{eq:elementary_interval}
		E_{\bm{\ell}, \bm{k}} = \prod_{j=1}^{s} \left[ \frac{{k}_j}{b^{{\ell}_j}},  \frac{{k}_j + 1}{b^{{\ell}_j}}\right),
	\end{equation}
	where $\bm{\ell} = ({\ell}_1, \dotsc, {\ell}_s) \in \mathbb{N}_0^s$, {$\abs{\bm{\ell}} := \sum_{j=1}^s {\ell}_j = m - t$} and $\bm{k} = ({k}_1, \dotsc, {k}_s) \in \mathbb{N}_0^s$ with ${k}_j < b^{{\ell}_j}$ for $j = 1, \dotsc, s$, contains exactly {$b^t$} points. Here, the smallest $t$ satisfying the above description is called the quality parameter, $s$ is the dimension, and $b^m$ is the number of quadrature points. Additional properties of $(t, m, s)$-nets can be found in~\cite{dick2010digital}. 

	To randomize digital nets, Owen proposes nested uniform scrambling~\cite{owen1995randomly} for the Sobol' sequence, although it requires substantial storage and computational costs. Matou\v sek~\cite{matouvsek1998thel2} proposes a random linear scramble as a more computationally efficient approach that does not alter the mean squared $L^2$ discrepancy. A review of other randomization methods can be found in~\cite{owen2003variance}. In this short note, we focus on the Owen-scrambled digital nets.

\subsection{Variance of the scrambled digital nets}
\label{sec:spectral_analysis_rqmc}
We denote by $I_N \coloneqq I_N (f)$ the Owen-scrambled $(t, m, s)$-net in base $b$ integration estimator for a general integrand $f$ defined on $[0, 1]^s$, where $N = b^{m}$. Following the derivations in~\cite{owen1997monte, dick2010digital}, the variance of the estimator ${I_N}$ is given by:
\begin{equation}
	\var{I_N} = \sum_{\bm{\ell} \in \mathbb{N}_0^s} \Gamma_{\bm{\ell}} \sigma^2_{\bm{\ell}},
\end{equation} 
where $\Gamma_{\bm{\ell}}$ is the gain coefficient, and $\sigma^2_{\bm{\ell}}$ represents the sum of squared Walsh coefficients over specific indices.

Next, we introduce the notation for the alternating sum. Given an interval $J = \prod_{j=1}^s [a_j, b_j]$, we define the alternating sum $\Delta (f, J)$ by
\begin{equation}
	\label{eq:alternating_sum}
	\Delta (f, J) = \sum_{\mathfrak{u} \subseteq 1:s} (-1)^{\abs{\mathfrak{u}}} f(\bm{a}^\mathfrak{u}: \bm{b}^{-\mathfrak{u}}),
\end{equation}
where $\bm{a} = (a_1, \dotsc, a_s)$ and $\bm{b} = (b_1, \dotsc, b_s)$. The expression $\bm{a}^\mathfrak{u}: \bm{b}^{-\mathfrak{u}}$ denotes the concatenation of two vectors such that the $j$-th component, $(\bm{a}^\mathfrak{u}: \bm{b}^{-\mathfrak{u}})_j$ equals $a_j$ if $j \in \mathfrak{u}$ and $b_j$ otherwise.

Using the alternating sum notation, we present the following definition.
\begin{Definition}[Generalized Vitali variation of order 2]
	In~\cite{dick2010digital}, the authors define the generalized Vitali variation of order $2$ as
	\begin{equation}
		\label{eq:generalized_variation_vitali}
		V_{\alpha} (f) = \sup_{\mathcal{P}} \left( \sum_{J \in \mathcal{P}} \mu(J) \abs*{\frac{{\Delta (f, J)}}{\mu(J)^{\alpha}} }^{2} \right)^{\frac{1}{2}},
	\end{equation}
	{where $0 < \alpha \leq 1$, the supremum is taken over all partitions $\mathcal{P}$ of $[0, 1]^s$ into axis-parallel subintervals, $\mu(J)$ denotes the Lebesgue measure of $J$, and $\Delta (f, J)$ is defined in~\eqref{eq:alternating_sum}.}
\end{Definition}
Following the derivations in~\cite{dick2010digital}, the variance $\var{I_N}$ can be bounded in terms of the generalized Vitali variation, as stated in the following proposition.
\begin{Proposition}[Variance of the scrambled $(t,m,s)$-net in base $b$ estimator]
	Following~\cite{dick2010digital}, the variance of the scrambled digital nets estimator $I_N$ satisfies
	\begin{equation}
		\var{I_{N}} \leq C_{\alpha, b, s, t} V^2_{\alpha}(f) N^{-1-2\alpha} \log^{s-1} N,
	\end{equation}
where the constant $C_{\alpha, b, s, t} < +\infty$ depends on $\alpha, b, s$ and $t$. 
\end{Proposition}
Let $\norm{\cdot}_p$ denote the $L^p$ norm on $[0, 1]^s$. Our goal is to show
\begin{equation*}
	V_{\alpha}(f) \leq \norm{\partial^{1:s} f}_{p},
\end{equation*}
where $\partial^{1:s} f$ denotes the weak mixed first-order derivative of $f$, $p = \frac{2}{3-2\alpha}$ for $\frac{1}{2} \leq \alpha \leq 1$. Consequently, the variance $\var{I_N}$ can be connected to the $L^p$ integrability of the weak derivative. For a general reference on weak derivatives, see~\cite{evans2022partial}. 

We present the following lemma connecting the alternating sum and the weak derivative.
\begin{Lemma}[Alternating sum and the weak derivative]
\label{lemma:weak_derivative_alternating_sum}
For a continuous integrand $f$ whose weak derivative $\partial^{1:s}f$ exists on an axis-parallel interval $J = \prod_{j=1}^s [a_j, b_j]$ with $0 \leq a_j < b_j \leq 1$, for $j = 1, \dotsc, s$, we have
\begin{equation}
	{{\Delta (f, J)}} = \int_{J} {\partial^{1:s} f} (\bm{t}) d\bm{t},
\end{equation}
where $\partial^{1:s} f$ denotes the weak derivative.
\end{Lemma}
Lemma~\ref{lemma:weak_derivative_alternating_sum} extends the results presented in~\cite{owen2005multidimensional}, where the derivative is defined in the usual sense on $J$. We present the proof below.
\begin{proof}
	When $f$ is continuous and weakly differentiable on $J$, Theorem 8.2 in~\cite{brezis2011functional} presents the following result for the 1-d case:
	\begin{equation}
		f(b) - f(a) = \int_{a}^{b} f^{\prime} (t) dt.
	\end{equation}
	We now proceed by induction. Decompose $J$ as $J=[a_1, b_1] \times J_{-1}$ with $J_{-1} = \prod_{j = 2}^s [a_j, b_j]$. Denote by $f(\cdot \mid t_1 = \tau)$ the restriction of $f(\bm{t})$ to $f(\tau, t_2, \dotsc, t_{s})$ for $\tau \in [0, 1]$. We have
	\begin{equation}
		\begin{split}
			{{\Delta (f, J)}} &= {\Delta (f(\cdot \mid t_1 = b_1), J_{-1})} - {\Delta (f(\cdot \mid t_1 = a_1), J_{-1})}\\
			&= \int_{a_1}^{b_1} \frac{\partial}{\partial \tau} \Delta (f(\cdot \mid t_1 = \tau), J_{-1}) d\tau\\
			&= \int_{a_1}^{b_1} \frac{\partial}{\partial \tau} \int_{J_{-1}} \frac{\partial f(\cdot \mid t_1 = \tau) }{\partial \bm{t}_{2:s}} d\bm{t}_{2:s} d \tau\\
			&= \int_{J} \partial^{1:s} f(\bm{t}) d\bm{t},
		\end{split}
	\end{equation}
	where in the third line we apply the induction hypothesis in $(s-1)$ dimensions, and in the fourth line we exchange the weak derivative and integration~\cite{cheng2006differentiation}. This concludes the proof. 
\end{proof}
In the following, we show that $V_{\alpha}(f)$ is bounded by $\norm{\partial^{1:s} f}_{p}$ with $p = \frac{2}{3-2\alpha}$. First, we apply H\"older's inequality to obtain
\begin{equation}
	{{\Delta (f, J)}} = \int_{J} \partial^{1:s}f(\bm{t}) d\bm{t} \leq  \norm{ \partial^{1:s} f \cdot \mathbbm{1}_{J}}_p \ \mu(J)^{1-\frac{1}{p}},
\end{equation}
where $\mathbbm{1}_{J}$ denotes the indicator function over the set $J$, taking the value 1 inside $J$ and 0 otherwise. Thus, for all partitions $\mathcal{P}$ of $[0, 1]^s$ into axis-parallel subintervals, we have
\begin{equation*}
	\sum_{J \in \mathcal{P}} \mu(J) \abs*{\frac{{\Delta (f, J)}}{\mu(J)^{\alpha}} }^{2} \leq \sum_{J \in \mathcal{P}} \mu(J)^{3-2\alpha - \frac{2}{p}} \norm{ \partial^{1:s} f \cdot \mathbbm{1}_{J}}_p^2.
\end{equation*}
Observe the decomposition of the $p$-th power of the $L^p$ norm over the partition $\mathcal{P}$:
\begin{equation*}
	\norm{ \partial^{1:s} f  }_p^p = \sum_{J \in \mathcal{P}} \norm{ \partial^{1:s} f \cdot \mathbbm{1}_{J}}_p^p.  
\end{equation*}
When $p \leq 2$, we have the following superadditivity condition:
\begin{equation}
		\left(\norm{ \partial^{1:s} f  }_p^p \right)^{\frac{2}{p}} = \left( \sum_{J \in \mathcal{P}} \norm{ \partial^{1:s} f \cdot \mathbbm{1}_{J}}_p^p \right)^{\frac{2}{p}} \geq  \sum_{J \in \mathcal{P}} \norm{ \partial^{1:s} f \cdot \mathbbm{1}_{J}}_p^2.
\end{equation}
Finally, we obtain
\begin{equation}
	V_{\alpha} (f) = \sup_{\mathcal{P}} \left( \sum_{J \in \mathcal{P}} \mu(J) \abs*{\frac{{\Delta (f, J)}}{\mu(J)^{\alpha}} }^{2} \right)^{\frac{1}{2}} \leq \sup_{\mathcal{P}} \left( \sum_{J \in \mathcal{P}} \norm{ \partial^{1:s} f \cdot \mathbbm{1}_{J}}_p^2 \right)^{\frac{1}{2}} \leq \norm{ \partial^{1:s} f  }_p.
\end{equation}
Thus, we establish a connection between the generalized Vitali variation and the $L^p$ integrability of the mixed first-order derivatives. Moreover, the convergence rate of the estimator can inferred from the $L^p$ the integrability of the weak derivative $\partial^{1:s} f$, for $1 \leq p \leq 2$. Since we consider the weak derivative, the analysis accommodates a broader class of continuous functions, such as the functions with kinks. Meanwhile, discontinuous functions are studied in~\cite{he2015convergence} and more recently in~\cite{liu2024randomized}. We conclude this short note by summarizing our contributions in the following remark.
\begin{remark}[On the convergence rate of the scrambled digital nets]
	Notice that while the generalized Vitali variation of order 2 is upper bounded by the $L^p$ norm of the derivative, the looser upper bound derived in this work provides a more tractable approach for determining the convergence rates. Specifically, our work extends prior studies in the following two aspects:
	\begin{itemize}
		\item The estimator variance convergence rate $\mathcal{O}(N^{-3} \log^{s-1} N)$ under weaker regularity: A sufficient condition is the $L^2$-integrability of the mixed first-order weak $\partial^{1:s} f$, which generalizes a sufficient condition derived in~\cite{dick2010digital} that requires the continuous derivative in the strong sense to be in $L^2$.
		\item Generalized conditions for rates $\mathcal{O}(N^{-2-\delta} \log^{s-1} N) (0 < \delta < 1)$: This work generalizes the boundary growth conditions considered in~\cite{liu2024randomized}, which model the behavior of the derivative near the boundary, to the integrability conditions on the weak derivatives.
	\end{itemize}
	These extensions can be useful for practitioners to determine the convergence rates of the scrambled digital nets for a broader class of integrands.
\end{remark}


\appendix

\section{Numerical Examples}
\label{sec:numex}
We present some numerical examples in the appendix. Specifically, we consider two integrands with kinks. All the numerical simulations use the Sobol' sequence with Matou\v{s}ek-scrambling~\cite{matouvsek1998thel2}, as implemented in the \texttt{scipy.qmc} module~\cite{2020SciPy-NMeth}. 
\subsection*{Example 1}
In Example 1, we consider an integrand of the form:
\begin{equation}
	f(\bm{t}) = \prod_{j=1}^s \abs*{t_j - \dfrac{1}{2}}^{\alpha}, \quad \alpha > 0,
\end{equation}
where the exact integration value is $\frac{1}{2^{\alpha s} (\alpha + 1)^s}$.
Notice that in this example, the integrand kinks are axis-parallel. When $0 < \alpha \leq \frac{1}{2}$, the weak derivative $\partial^{1:s} f$ is in $L^{-\frac{1}{\alpha - 1} - \epsilon}$ for any arbitrarily small $\epsilon > 0$ and the convergence rate for the estimator variance is $\mathcal{O}(N^{-2-{2}{\alpha} + \epsilon} \log^{s-1} N)$ for any arbitrarily small $\epsilon > 0$.
\begin{figure}[htbp]
	\centering
	\begin{subfigure}{0.32\textwidth}
		\includegraphics[width=\textwidth]{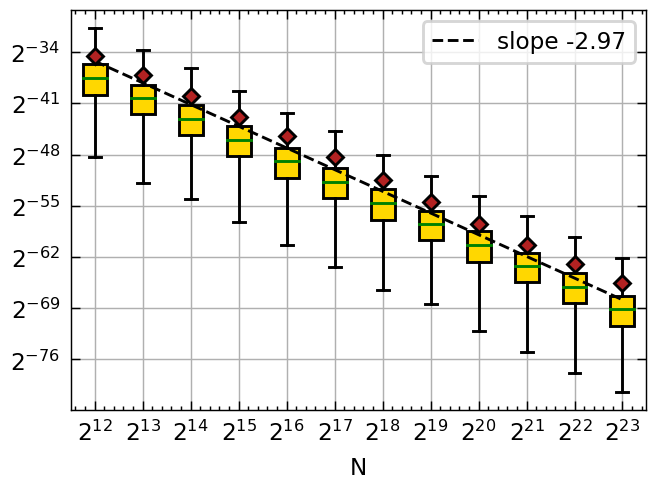}
		\caption{$s = 2, \alpha = \frac{1}{2}$.}
	\end{subfigure}
	\begin{subfigure}{0.32\textwidth}
		\includegraphics[width=\textwidth]{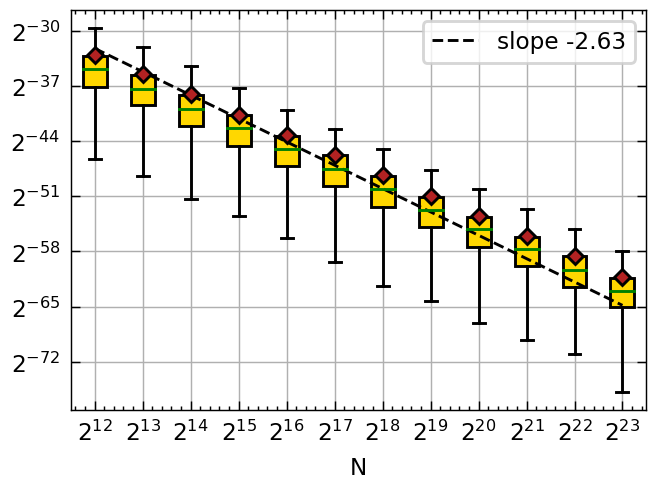}
		\caption{$s = 2, \alpha = \frac{1}{3}$.}
	\end{subfigure}
	\begin{subfigure}{0.32\textwidth}
		\includegraphics[width=\textwidth]{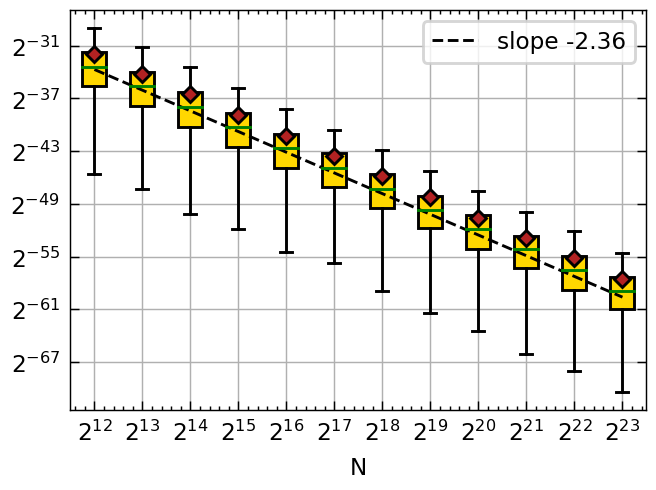}
		\caption{$s = 2, \alpha = \frac{1}{5}$.}
	\end{subfigure}
	\\
	\begin{subfigure}{0.32\textwidth}
		\includegraphics[width=\textwidth]{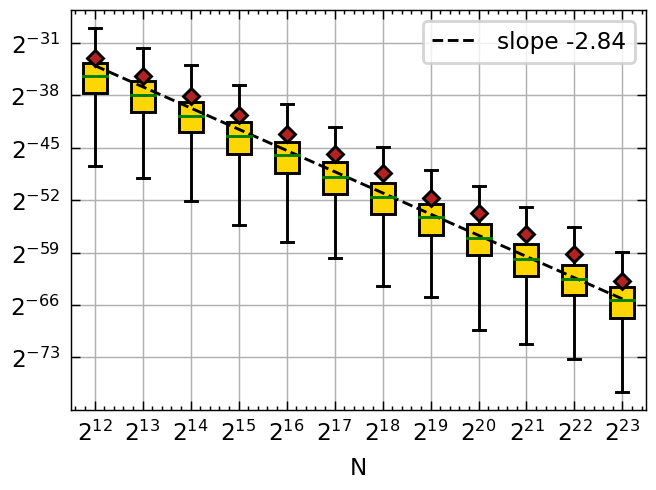}
		\caption{$s = 3, \alpha = \frac{1}{2}$.}
	\end{subfigure}
	\begin{subfigure}{0.32\textwidth}
		\includegraphics[width=\textwidth]{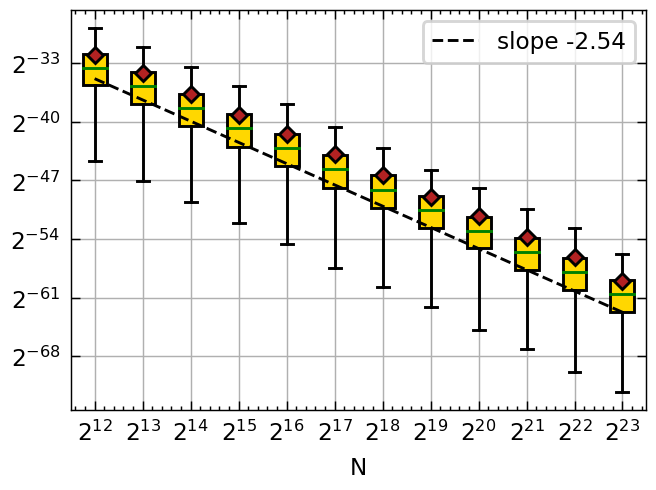}
		\caption{$s = 3, \alpha = \frac{1}{3}$.}
	\end{subfigure}
	\begin{subfigure}{0.32\textwidth}
		\includegraphics[width=\textwidth]{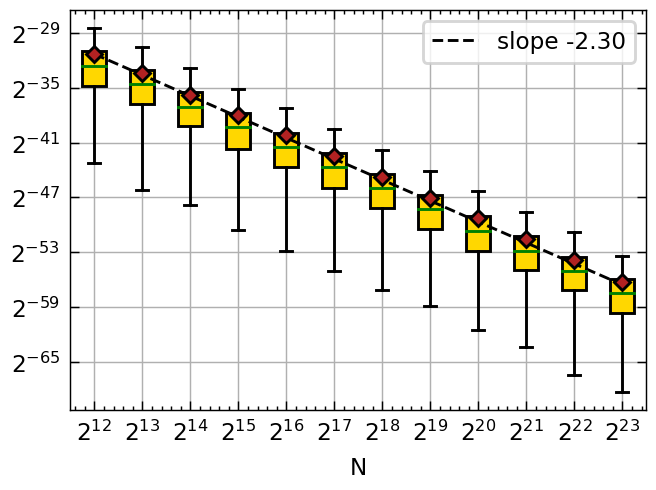}
		\caption{$s = 3, \alpha = \frac{1}{5}$.}
	\end{subfigure}
	\\
	\begin{subfigure}{0.32\textwidth}
		\includegraphics[width=\textwidth]{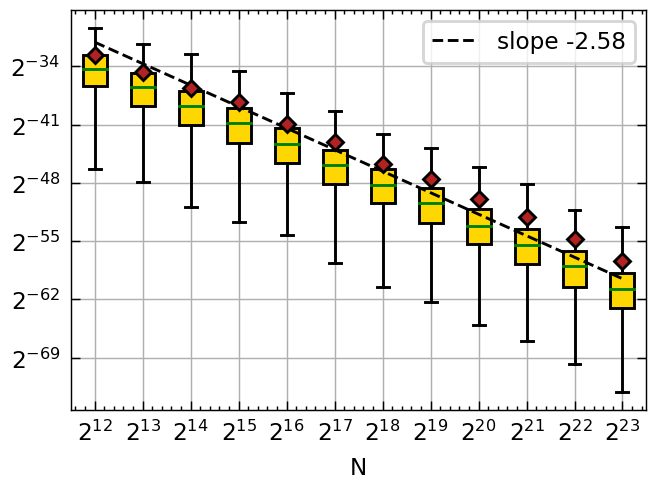}
		\caption{$s = 5, \alpha = \frac{1}{2}$.}
	\end{subfigure}
	\begin{subfigure}{0.32\textwidth}
		\includegraphics[width=\textwidth]{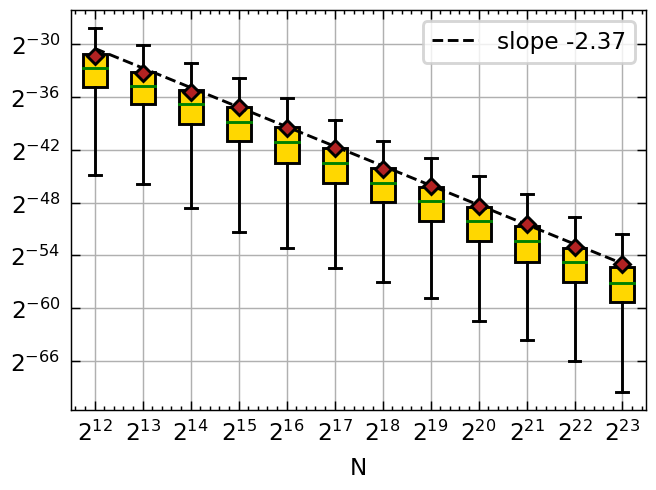}
		\caption{$s = 5, \alpha = \frac{1}{3}$.}
	\end{subfigure}
	\begin{subfigure}{0.32\textwidth}
		\includegraphics[width=\textwidth]{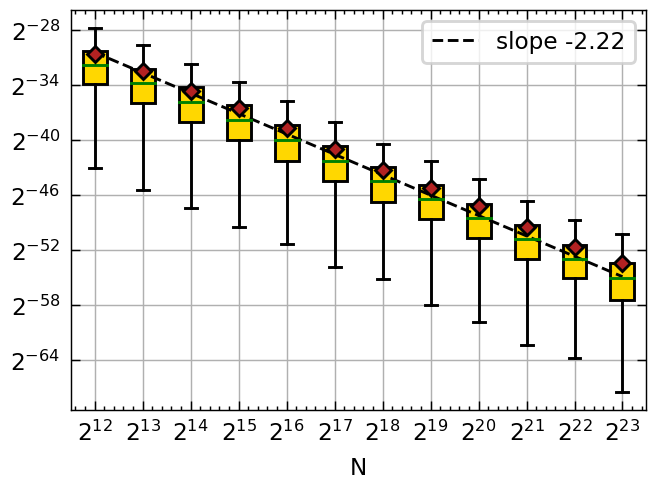}
		\caption{$s = 5, \alpha = \frac{1}{5}$.}
	\end{subfigure}
	\caption{Example 1: The boxplot characterization of squared error distributions, $(I_N - I)^2$, for scrambled Sobol' sequence estimators across various dimensions $s$ and parameters $\alpha$. Each whisker in the boxplot extends from the 1st to 99th percentile of 8,192 independent realizations of the squared errors.}
	\label{fig:numerical_results_example1}
\end{figure}

Figure~\ref{fig:numerical_results_example1} presents the squared errors $(I - I_N)^2$ of the scrambled digital net estimators for a range of $N$ for various dimensions $s$ and parameter $\alpha$. The empirical convergence rates for the cases $s=2$ with $\alpha = \frac{1}{2}, \frac{1}{3}, \frac{1}{5}$ are close to the asymptotic rates $\mathcal{O}(N^{-2-2\alpha + \epsilon})$, where one can almost neglect the effect from the logarithmic in the complexity. For the cases $s=3$ and $s=5$, there are more nonasymptotic effects from the logarithmic term in the complexity.
\subsection*{Example 2}
In Example 2, we consider an integrand of the form:
\begin{equation}
	f(\bm{t}) = \max\left(\sum_{j=1}^s t_j - 1, 0\right)^{s + \alpha}, \quad \alpha > -1.
\end{equation}
In this example, the kinks are along the hyperplane $\sum_{j=1}^s t_j = 1$. When $-1 < \alpha \leq -\frac{1}{2}$, the weak derivative $\partial^{1:s}$ is in $L^{-\frac{1}{\alpha} - \epsilon}$ for any arbitrarily small $\epsilon > 0$ and the convergence rate for the estimator variance is $\mathcal{O}(N^{-4-{2}{\alpha} + \epsilon} \log^{s-1} N)$ for any arbitrarily small $\epsilon > 0$.

\begin{figure}[htbp]
	\centering
	\begin{subfigure}{0.32\textwidth}
		\includegraphics[width=\textwidth]{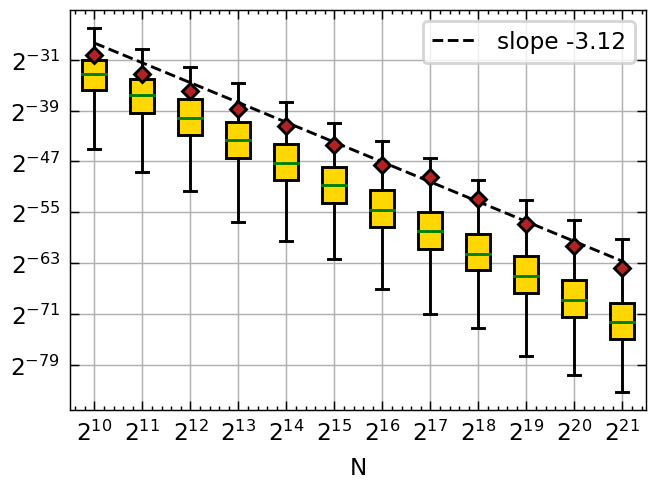}
		\caption{$s = 2, \alpha = -0.5$.}
	\end{subfigure}
	\begin{subfigure}{0.32\textwidth}
		\includegraphics[width=\textwidth]{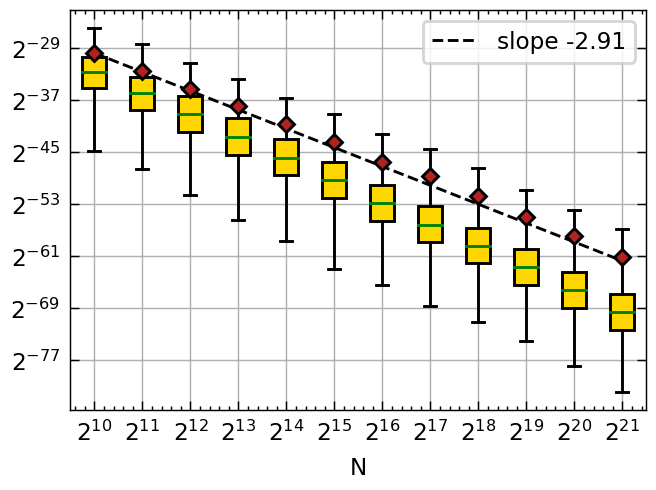}
		\caption{$s = 2, \alpha = -0.7$.}
	\end{subfigure}
	\begin{subfigure}{0.32\textwidth}
		\includegraphics[width=\textwidth]{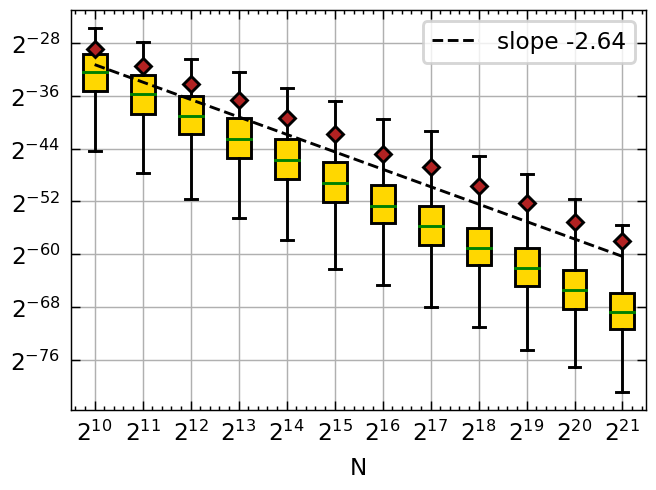}
		\caption{$s = 2, \alpha = -0.9$.}
	\end{subfigure}
	\\
	\begin{subfigure}{0.32\textwidth}
		\includegraphics[width=\textwidth]{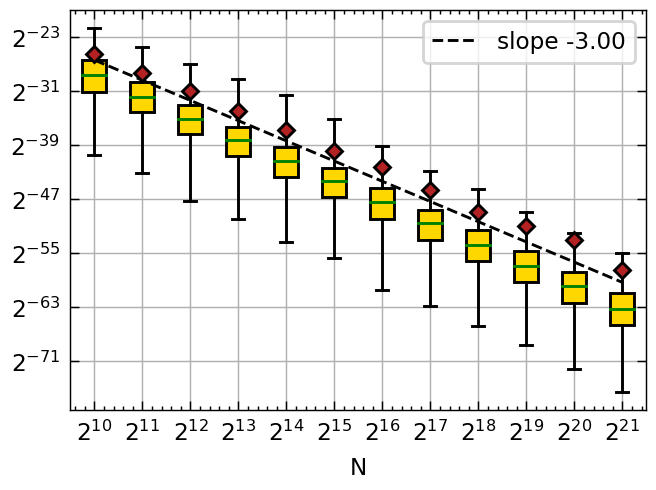}
		\caption{$s = 3, \alpha = -0.5$.}
	\end{subfigure}
	\begin{subfigure}{0.32\textwidth}
		\includegraphics[width=\textwidth]{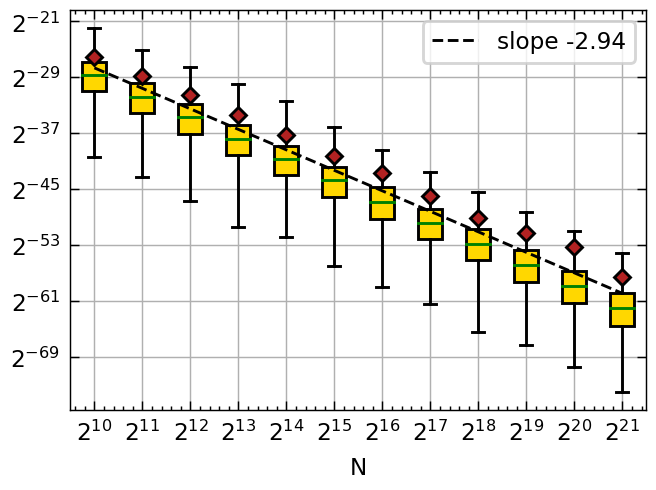}
		\caption{$s = 3, \alpha = -0.7$.}
	\end{subfigure}
	\begin{subfigure}{0.32\textwidth}
		\includegraphics[width=\textwidth]{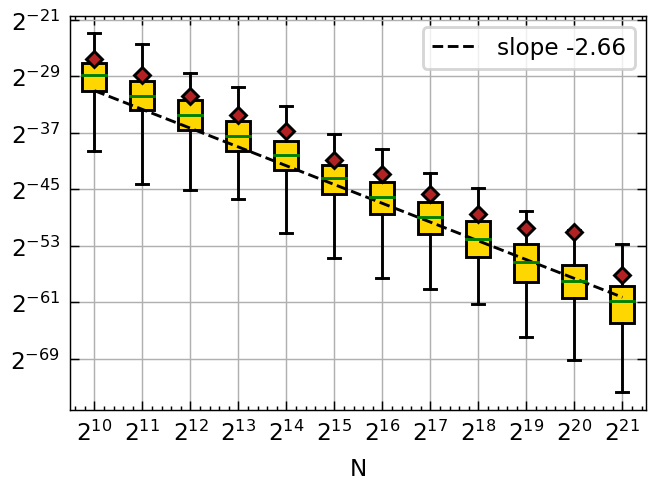}
		\caption{$s = 3, \alpha = -0.9$.}
	\end{subfigure}
	\\
	\begin{subfigure}{0.32\textwidth}
		\includegraphics[width=\textwidth]{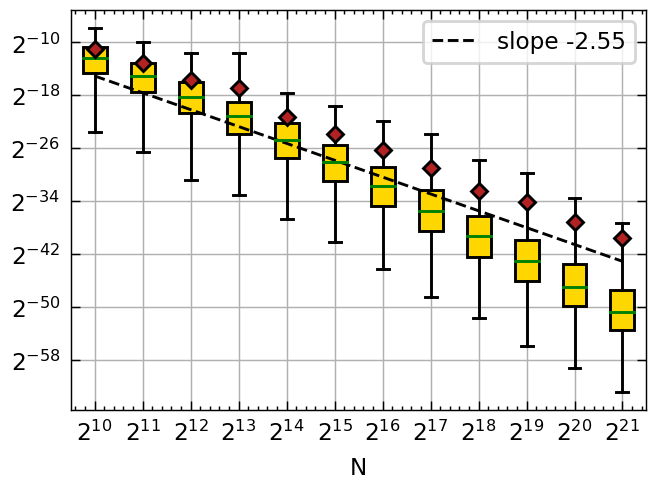}
		\caption{$s = 5, \alpha = -0.5$.}
	\end{subfigure}
	\begin{subfigure}{0.32\textwidth}
		\includegraphics[width=\textwidth]{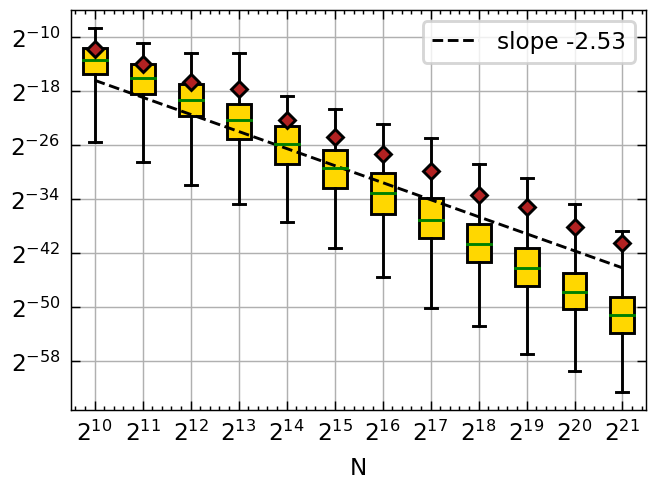}
		\caption{$s = 5, \alpha = -0.7$.}
	\end{subfigure}
	\begin{subfigure}{0.32\textwidth}
		\includegraphics[width=\textwidth]{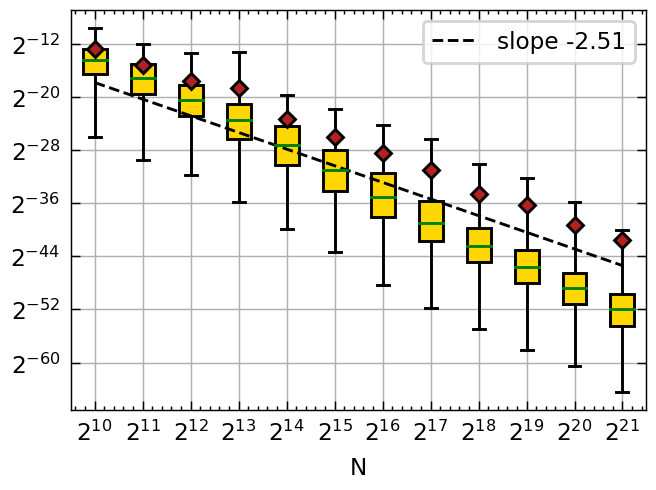}
		\caption{$s = 5, \alpha = -0.9$.}
	\end{subfigure}
	\caption{Example 2: The boxplot characterization of squared error distributions, $(I_N - I)^2$, for scrambled Sobol' sequence estimators across various dimensions $s$ and parameters $\alpha$. The reference value $I$ is approximated by averaging 8,192 independent realizations of scrambled Sobol' sequence estimators with quadrature size $N = 2^{25}$. Each whisker in the boxplot extends from the 1st to 99th percentile of 8,192 independent realizations.}
	\label{fig:numerical_results_example2}
\end{figure}

Figure~\ref{fig:numerical_results_example2} presents the numerical results for Example 2, where the reference values are computed with an ensemble average of 8,192 independent realizations of scrambled Sobol' sequence estimator with quadratures size $N = 2^{25}$. For the cases $s = 2$ and $s = 3$, when $\alpha = -0.5$, empirical rates approach $\mathcal{O}(N^{-3 + \epsilon})$, aligning theoretical predictions. For the cases $\alpha = -0.7$ and $\alpha = -0.9$, the convergence rates exceed our a priori estimate rates, which suggests the upper bound derived in this work may not be tight. More detailed analysis for this type of integrand is left for future work. For the case $s = 5$, the nonasymptotic effects become more profound due to the increased integration dimension.

\bibliographystyle{siam}
\bibliography{sampling, bibliography_QMC_theory, bibliography_QMC_finance, hoqmc}


\end{document}

%% file: commands_qmc.tex
\usepackage{mathtools}

\providecommand{\var}[1]{{\ensuremath{\mathrm{Var}}\mspace{-2mu}\left[#1\right]}}

\DeclarePairedDelimiter\abs{\lvert}{\rvert}%
\DeclarePairedDelimiter\norm{\lVert}{\rVert}%

\newtheorem{Definition}{Definition}

\newtheorem{Proposition}{Proposition}
\newtheorem{Lemma}{Lemma}

\newcommand{\RN}[1]{%
	\textup{\uppercase\expandafter{\romannumeral#1}}%
}